\documentclass[a4paper]{article}
\usepackage[colorlinks=true,citecolor=black,linkcolor=black]{hyperref}
\usepackage{geometry}
\usepackage{amsmath,amssymb,amsthm}
\usepackage{color}
\usepackage{hyperref, doi}
\usepackage{subfig} %
\usepackage{graphicx}
\usepackage{multirow}

\usepackage{booktabs} 
\usepackage{bm}
\usepackage{array} 
\usepackage{paralist} 
\usepackage{verbatim} 
\usepackage[font=small]{caption}

\newcommand{\ZZ}{\mathbf{Z}}
\newcommand{\mP}{\mathcal{P}}

\DeclareMathOperator{\Div}{div}

\DeclareMathOperator{\argmin}{argmin}

\newcommand{\grad}{\nabla}

\numberwithin{equation}{section}

\newtheorem{proposition}{Proposition}[section]

\theoremstyle{remark}
\newtheorem{remark}{Remark}[section]
\theoremstyle{definition}

\newcommand{\bitem}{\begin{itemize}}
	\newcommand{\eitem}{\end{itemize}}

\newcommand{\bpm}{\begin{pmatrix}}
	\newcommand{\epm}{\end{pmatrix}}

\newcommand{\bq}{\begin{equation}}
\newcommand{\eq}{\end{equation}}

\newcommand{\D}{\mathrm{div} \;}

\let\abs=\envert

\let\norm=\enVert

\let\inprod=\inProd

\begin{document}
\title{Variational Image Motion Estimation by Accelerated Dual Optimization}
\author{Hongpeng Sun\thanks{Institute for Mathematical Sciences,
Renmin University of China, No.~59, Zhongguancun Street, Haidian District,
100872 Beijing, People's Republic of China.
Email: \href{mailto:hpsun@amss.ac.cn}{hpsun@amss.ac.cn}.} \quad 
 Xuecheng Tai\thanks{Department of Mathematics, Hong Kong Baptist University,
			Kowloon Tong, Hong Kong.
			Email: \href{mailto:xuechengtai@hkbu.edu.hk}{xuechengtai@hkbu.edu.hk}.}
 \quad
	Jing Yuan \thanks
	{School of Mathematics and Statistics, Xidian University,  
		Email: \href{mailto:jyuan@xidian.edu.cn}{jyuan@xidian.edu.cn}.}
}


\maketitle

\begin{abstract}
    Estimating optical flows is one of the most interesting problems in computer vision, which estimates the essential information about pixel-wise displacements between two consecutive images. This work introduces an efficient dual optimization framework with accelerated preconditioners to the challenging nonsmooth optimization problem of total-variation regularized optical-flow estimation. In theory, the proposed dual optimization framework brings an elegant variational analysis on the given difficult optimization problem, while presenting an efficient algorithmic scheme without directly tackling the corresponding nonsmoothness in numeric.
	By introducing efficient preconditioners with a multi-scale implementation, 
	the proposed accelerated dual optimization approaches achieve 
	competitive estimation results of image motion, comparing to the state-of-the-art methods. 
	Moreover, we show that the proposed preconditioners can guarantee convergence of the implemented numerical schemes with high efficiency.

\end{abstract}


%

         
\section{Introduction}
\label{intro}
Accurate estimation of optical flows is one of the most interesting but still challenging problems in both mathematical image processing and computer vision \cite{BPS,FW,FBK, RS, ADK}, which aims to discover the essential information about pixel-wise displacements $\bm{d}:=(u,v)^T$ within the given two consecutive images $I(x,y, t)$ and $I(x,y, t+1)$. The basic idea that the intensity of a point keeps constant along its trajectory \cite{ADK}, where the intensity values of the corresponding pixels in the two given images are the same, is the so-called \emph{optical-flow constraint}:
\[
I(x+u, y + v, t+1) - I(x,y, t) \, \simeq \, \underbrace{\partial_x I\, u  + \partial_y I\, v}_{\nabla I^{T} \bm{d}} + \partial_t I \, = \, 0 \,
\]
given that the displacement $\bm{d}$ at each pixel is small enough.

Energy minimization approach based on variational functional is one of the most important methods to tackle this problem \cite{ADK, BA,BR, SRB, ZPH}.
We mainly focus on the variational functional based method, which is based on the following  common starting point: the \emph{optical flow constraint},
\begin{equation}
\frac{dI(x,y,t)}{dt} = 0,
\end{equation}
where $I(x,y,t)$ denotes the image intensity (or brightness)  at the
point $(x, y)$ in the image plane at time $t$ \cite{HS}. The the \emph{optical flow constraint} means that the intensity of a point keeps constant along its trajectory \cite{ADK}. If the displacements of $I(x,y,t)$ are very small, with the chain rule for differentiation on $t$, we can get the first order approximation of the \emph{optical flow constraint},
\begin{equation}\label{eq:ofc:diff}
\frac{\partial I}{\partial x} \frac{\partial x}{\partial t} + \frac{\partial I}{\partial y} \frac{\partial y}{\partial t} + \frac{\partial I}{\partial t}=0.
\end{equation}
Denoting $\nabla =[\nabla_x, \nabla_y]^T$, $\nabla I = (\frac{\partial I}{\partial x}, \frac{\partial I}{\partial y})^T:  =[I_x,I_y]^T$, $u := \frac{\partial x}{\partial t}$, $v:=\frac{\partial y}{\partial t}$,  $\bm{d}:=(u,v)^T$ and $I_t: = \frac{\partial I}{\partial t}$. With additional smooth quadratic regularization of $\nabla u$, $\nabla v$, the classical Horn--Schunck functional is as follows
\begin{equation}
\int_{\Omega}(\nabla I^T\bm d+I_t)^2d\sigma + \lambda( |\nabla u|^2 + |\nabla v|^2)d\sigma,
\end{equation}
where $\Omega$ is the image domain, $\lambda >0$ is a regularization parameter, $d\sigma=dxdy$ is the area element of $\Omega$, and $\bm d = (u,v)^T$ is the optical flow vector. It can  yield high density of flow vectors and is still very useful for lots of applications \cite{SRB}. However, due to the well-known drawback of the $L^2$ quadratic data or regularization terms, it does not take into account the discontinuities of flow field and thus does not  allow for discontinuities in the displacement field. It also  does not allow for outliers in the data term and is very sensitive to noise. It is one of the earliest variational functional for optical flow. A lot of variational frameworks were developed since then. One is the  robust estimation framework \cite{BA, BR}. By introducing the robust estimators from the robust statistics, the following variational framework is built in \cite{BA, BR},
\begin{equation}
\int_{\Omega}\rho_{D}(\nabla I^T\bm d+I_t)d\sigma + \lambda( \rho_{S}(|\nabla_x u|)+ \rho_{S}(|\nabla_y u|) + \rho_{S}(|\nabla_x v|) + \rho_{S}(|\nabla_y v|))d\sigma,
\end{equation}
where $\rho_{D}(\cdot)$ and $\rho_{S}(\cdot)$ are robust estimators which are usually nonconvex. Quite a lot of the robust estimators $\rho_{S}(\cdot)$  have the corresponding Markov random field interpretations \cite{BA}. The convergence analysis of some iteration algorithms for this framework in infinite dimensional spaces can be found in \cite{ADK}.
Variational nonlocal model is also very useful for optical flow estimates \cite{KKO,RBP, WPB}.
Another widely employed framework is the following convex TV-$L^1$ variational framework, which is our focus \cite{ZPH,WPZBC,PMF},
\begin{equation}\label{eq:l1:optical:flow:zach}
\min_{\bm{d}}  \int_{\Omega} ( |Du|+|Dv|) + \lambda \int_{\Omega} |\rho(\bm d)| d\sigma,
\end{equation} 
where $\rho(x)$ is defined as
\[
\rho(\bm d) :=\nabla I^T\bm d+I_t , \quad \bm d = (u,v)^{T}.
\]

For the algorithmic development, alternating minimization method is developed for \eqref{eq:l1:optical:flow:zach} with adding the penalty terms, i.e., 
\begin{equation}\label{eq:alter:appro}
E_{\theta}(\bm d, \bm h):  = \int_{\Omega} ( |Du|+|Dv|) +  \lambda \int_{\Omega} |\rho(\bm h)| d\sigma + \int_{\Omega}\frac{|\bm d-\bm h|^2}{2 \theta}d\sigma,
\end{equation}
which can be seen as a penalty method for the following constrained optimization problem without increasing ${1}/{\theta}$ during iterations \cite{NW},
\begin{equation}\label{eq:l1:optical:flow:zach:penallty}
\min_{\bm{d}, \bm h}  \int_{\Omega} ( |Du|+|Dv|) + \lambda \int_{\Omega} |\rho(\bm h)| d\sigma, \quad \bm d = \bm h.
\end{equation} 
The model \eqref{eq:alter:appro} can be seen as an approximation to \eqref{eq:l1:optical:flow:zach} with fixed $\theta$. The choice of $\theta$ thus becomes a subtle issue \cite{PMF}.
Given $\bm d^0$ or $\bm h^0$, one can minimize $\bm d$ and $\bm h$ alternatively in \eqref{eq:alter:appro}. With fixed $\bm d^k$, the minimization problem of $\bm h$ has explicit solution; with fixed $\bm h^{k+1}$, the minimization problem of $\bm d$ turns out to be the usual ROF denoising problem, where various optimization methods can be employed \cite{ZPH, PMF}.

In this work, we mainly focus on a new isotropic TV-$L^1$ optical flow model, i.e., the following illumination-invariant TV-$L^1$ optical flow model \cite{CP}. This optical flow model introduces an additional function $w$, standing for varying illumination, into the optical flow constraint such that
\begin{equation}\label{eq:l1:optical:flow}
\min_{w,\bm d} \; \int_{\Omega} |Dw| \, + \, \int_{\Omega}  |D \bm d| \, + \, \lambda \int_{\Omega}|\rho(w,\bm d)|d\sigma,
\end{equation} 
where $\lambda$ is a positive parameter, and the new illumination invariant optical flow constraint $\rho(w,x)$ is defined as
\[
\rho(w,\bm d)\, := \, \nabla I^T\bm d+ \beta w +I_t \, , 
\]
$\beta$ is the positive parameter to balance the introduced illumination term $w$. Such illumination-invariant optical flow model \eqref{eq:l1:optical:flow} is motivated by the observation that image intensities do not strictly stay constant over time due to illumination changes and shadows in many real-world situations \cite{CP, SH}. 
Clearly, the following classical TV-$L^1$ isotropic optical flow model  can be recovered directly from \eqref{eq:l1:optical:flow} while setting $\beta=0$ which is similar to the anisotropic case as in \eqref{eq:l1:optical:flow:zach} \cite{ZPH,PMF, WPZBC}
\begin{equation}\label{eq:l1:optical:flow:no:w}
   \min_{\bm d} \; \, \int_{\Omega}  |D \bm d| \, + \, \lambda \int_{\Omega}|\rho(\bm d)|d\sigma, \quad \rho(\bm d)\, := \, \nabla I^T\bm d +I_t \, .
\end{equation}

With the recent developments of convex optimization algorithms, e.g.,  \cite{CP}, a lot of methods can be used to solve \eqref{eq:l1:optical:flow} and \eqref{eq:l1:optical:flow:no:w}. First-order primal-dual method is employed for \eqref{eq:l1:optical:flow} in \cite{CP}. There exist some ADMM algorithms for the optical flow estimate; see \cite{FSU} with piecewise-affine model and  \cite{ZK} with $l_p$ regularization. We also refer to \cite{VGSG1,VGSG2} for image registration with ADMM   which is quite similar to the optical flow estimate. ADMM was employed to primal problems in these work.  There are usually more than two primal variables and all primal variables are updated consecutively; see \cite{FSU}. However, it was recently found that ADMM can diverge even for convex problems if there are more than 2 block of variables  \cite{CBHY}. Another problem of the primal ADMM is that complicated nonlinear subproblems usually appeared for some subproblems, where gradient-based minimization techniques have to be used for dealing with the corresponding subproblems \cite{VGSG2}.

In this paper, we proposed a convergent and preconditioned ADMM for the optical flow estimates \eqref{eq:l1:optical:flow} and \eqref{eq:l1:optical:flow:no:w} via the dual approach. Our contributions are as follows. We first proposed a dual framework for the optical flow estimates  \eqref{eq:l1:optical:flow} and \eqref{eq:l1:optical:flow:no:w}.  For the dual problem  \eqref{eq:l1:optical:flow:no:w}, although there are only two block of variables, however, each subproblem is still highly nonlinear and is hard to solve.  We circumvent this problem through preconditioning techniques. 
For the dual problem of  \eqref{eq:l1:optical:flow}, there are three blocks of dual variables. We regroup the three block of variables into two big blocks. Special block preconditioners are designed for each block variables. Convergence can be guaranteed by the preconditioned ADMM  framework \cite{BS,SUN} and \cite{DY,LST,SUN}. Besides, we also studied the preconditioned Douglas-Rachford splitting method for \eqref{eq:l1:optical:flow} and \eqref{eq:l1:optical:flow:no:w} based on the corresponding primal-dual form. Furthermore, we studied the over-relaxed variants of these preconditioned ADMM. It turns out that these relaxation variants can bring out certain accelerations \cite{DY,FG, LST,SUN} for some case instead of all cases as  shown in \cite{DY, LST, SUN},   probably due to the multiscale pyramid structure.  
To the best knowledge of the authors, the preconditioned ADMM algorithms proposed in this paper are the first convergent variants of ADMM without error control for both \eqref{eq:l1:optical:flow} and \eqref{eq:l1:optical:flow:no:w} through the dual approach in literature.

The remaining of the paper is organized as follows. In section \ref{sec:primal-dual}, we give the primal-dual and dual models of the illumination-invariant TV-$L^1$ optical flow model, where the case \eqref{eq:l1:optical:flow:no:w} can also be covered.  In section \ref{sec:framework}, we give an introduction of the TV-$L^1$ variational framework and the preconditioned ADMM algorithms of optical flow estimates. With preconditioners or block preconditioners, we can get convergent and efficient ADMM for optical flow estimates via the dual approach. We also present a preconditioned Douglas-Rachford splitting method based on the saddle-point structure of the corresponding model.  In section \ref{sec:num}, we present a detailed numerical comparison with the proposed algorithms. In section \ref{sec:conclusion}, we give  some comments and a final conclusion. In section \ref{sec:appendix}, we give some additional details in the precious sections.

\section{Illumination-Invariant TV-$L^1$ Optical Flow Model}\label{sec:primal-dual}

\subsection{Equivalent Primal-Dual Model}

Now we study the new illumination-invariant optical flow model \eqref{eq:l1:optical:flow}, also called the \emph{primal model} in this work.  
For the total-variation regularization term $\int_{\Omega} |Dw|\,$ of \eqref{eq:l1:optical:flow}, we have its equivalent dual formulation \cite{CTV} such that
\bq \label{eq:tvw}
\int_{\Omega} |Dw| \, = \, \max_{p}\, \int_{\Omega} \nabla w \cdot p \, d\sigma \, = \, \max_{p}\,\int_{\Omega} -w \cdot \D p \, d\sigma \, , \quad \text{s.t.} \;\; \| p \|_{\infty} \leq 1 \, ,
\eq
where the dual variable $p(x):=(p_1(x),p_2(x))$ is a vector function in $\Omega$
with a vanishing boundary condition, i.e. $p|_{\partial \Omega} = 0$, and its infinity-norm $ \| p \|_{\infty}:=\max_{x\in\Omega} \sqrt{p_1^2(x) + p_2^2(x)}$ is less than $1$. Henceforth, $-\D:= \nabla^*$ is the adjoint operator of $\nabla$.

Likewise, we also have
\bq  \label{eq:tvu}
\int_{\Omega} |D \bm d| \, = \,\max_{q}\, \int_{\Omega} (\nabla u,\nabla v)^T \cdot (q_1,q_2)^T \, d\sigma \, = \, \max_{q_1}\,\int_{\Omega} -(u \cdot \D q_1 + v \cdot \D q_2) \, d\sigma \, ,
\eq
with
\bq \label{eq:tvv}
q(x) := (q_1(x), q_2(x)), \quad q_1(x) = (q_1^1(x), q_1^2(x)), \quad  q_2(x) = (q_2^1(x), q_2^2(x)),  \quad \text{s.t.} \;\; \| q \|_{\infty} \leq 1 \, , 
\eq
along with vanishing boundaries of the two vector functions $q_1(x)$ and $q_2(x)$ and the infinity-norm of $q$  defined as follows
\[
\| q \|_{\infty}:=\max_{x\in\Omega} \sqrt{(q_1^1)^2(x) + (q_1^2)^2(x) +(q_2^1)^2(x) + (q_2^2)^2(x)  }.
\]
Moreover, by the dual representation of the absolute function, the illumination-invariant optical flow fidelity term of \eqref{eq:l1:optical:flow} can be equally written as
\bq \label{eq:ofc-dual}
\lambda \int_{\Omega}|\rho(w,\bm d)|\, d\, \sigma 
\,=\, \max_s \,  \int_{\Omega}(\nabla I^T\bm d+ \beta w +I_t ) s\, d\sigma \, , \quad
\text{s.t.} \;\; \norm{s}_{\infty} \leq \lambda
\eq
where the infinity-norm $\norm{s}_{\infty} \leq \lambda$ means that $\abs{s(x)} \leq \lambda$ for $\forall x \in \Omega$.

In view of the equivalent formulations \eqref{eq:tvw}-\eqref{eq:ofc-dual}, the optimization problem \eqref{eq:l1:optical:flow} of the illumination-invariant optical-flow model can thus be identically expressed as
\bq \label{eq:pd1}
\min_{w, u, v} \max_{p, q_1,q_2,s}\, \int_{\Omega}(\nabla I^T\bm d+ \beta w +I_t ) s\, d\sigma \, - \, \int_{\Omega} w \cdot \D p \, d\sigma \, - \, \int_{\Omega} u \cdot \D q_1 \, d\sigma \, - \, \int_{\Omega} v \cdot \D q_2 \, d\sigma
\eq
subject to
\bq \label{eq:pd1-cond}
\| p \|_{\infty} \leq 1\, ,\;\;  \| q \|_{\infty} \leq 1\, ,\;  \| s \|_{\infty} \leq \lambda\, .
\eq
It is obvious that the minimax optimization model \eqref{eq:pd1} is equivalent to the primal optimization problem \eqref{eq:l1:optical:flow} by the Fenchel-Rockafellar duality theory (see Chapter 4.3 of \cite{KK}), which is called the \emph{primal-dual} model to \eqref{eq:l1:optical:flow} and can be further re-organized as
\bq \label{eq:pd2}
\min_{w, u, v} \max_{p, q_{1},q_2,s}\, \int_{\Omega} I_t \cdot s\, d\sigma\, + \, \int_{\Omega} w \cdot (\beta s - \D p) \, d\sigma \, + \, \int_{\Omega} u \cdot (sI_x - \D q_1) \, d\sigma \, + \, \int_{\Omega} v \cdot (sI_y - \D q_2) \, d\sigma
\eq
subject to the constraints \eqref{eq:pd1-cond} on the dual variables ($s, p, q_1, q_2$).

\subsection{Equivalent Dual Model}

Now we discretize and analyze the optimization problem \eqref{eq:l1:optical:flow} of optical-flow and its mathematically identical primal-dual representation \eqref{eq:pd2} in finite dimensional spaces. 
After discretization, we denote $X$ as the finite-dimensional space for any scalar image function $f_{i,j}$ where $(i,j) \in \Omega$, along
with the standard inner product for two scalar functions  $f, g \in X=\mathbb{R}^{M\times N}$ such that
\[
\inprod{f, g}_X \, = \, \sum_{(i,j) \in  \Omega} f_{i,j}  g_{i,j} \, ;
\]
also, the finite-dimensional vector space $Y = X \times X$. Here $\Omega \in \ZZ^2 $ denotes the following discretized grid \cite{CP}
\[
\Omega := \{(i,j)\ |\ i,j \in \mathbb{N}, \ 1 \leq i \leq M,
1 \leq j \leq N\}. 
\]
For \eqref{eq:l1:optical:flow}, we arrive at
\begin{equation}\label{eq:l1:finite}
\min_{w \in W} \min_{u,v \in  X} \|\nabla w\|_{1} + \|\nabla \bm d  \|_{1} +  \lambda \|\rho(w,\bm d)\|_{1},
\end{equation}
where the finite dimensional Hilbert spaces $W$ and $X$ are the corresponding image spaces.
The discrete  gradient operator $\grad: X \to Y$  can be found in \cite{CP, BS}.  Denoting $Y = X\times X$, with the standard scalar product, the discrete  divergence is  the negative adjoint of $\nabla$ with finite difference method \cite{CP, BS},
i.e., the unique linear mapping $\Div: Y \to X$ which satisfies
\[
\langle \nabla u, p\rangle_{Y} = \langle u, \nabla^{*} p\rangle_{X} = -\langle u,   \Div p \rangle_{X}, \quad \forall u \in X, \ p \in Y.
\]
Similarly, for the discrete version of \eqref{eq:l1:optical:flow:no:w}, we have
\begin{equation}\label{eq:l1:finite:no:w}
\min_{ \bm d \in  Y}  \|\nabla \bm d  \|_{1} +  \lambda \|\rho(\bm d)\|_{1}.
\end{equation}
We see that, in the primal-dual model \eqref{eq:pd2}, the variables ($w,u,v$) are free and their variations give rise to the following maximization problem:
\bq \label{eq:dual}
\max_{p, q_{1},q_2,s}\, \int_{\Omega} I_t \cdot s\, d\sigma\, 
\eq
subject to
\bq \label{eq:dual-cond}
\beta s - \D p \, =\, 0\, , \quad s I_x - \D q_1\,=\, 0\,, \quad sI_y - \D q_2 = 0 \, ,
\eq
along with the constraints as in \eqref{eq:pd1-cond}. Clearly, the optimization model \eqref{eq:dual} is also equivalent to the primal optimization problem \eqref{eq:l1:optical:flow} and its corresponding primal-dual model \eqref{eq:pd1} and \eqref{eq:pd2}. We thus called the model \eqref{eq:dual} as the \emph{dual model} in this study.

\section{TV-$L^1$ optical flow and preconditioned ADMM}\label{sec:framework}
In this section, we will focus on the discrete TV-$L^1$ models \eqref{eq:l1:finite}, \eqref{eq:l1:finite:no:w} and the corresponding preconditioned ADMM methods. We will also discuss the preconditioned Douglas-Rachford splitting method which is also very efficient. Let's begin with the primal-dual and dual reformulations of   \eqref{eq:l1:finite} and \eqref{eq:l1:finite:no:w}.
\subsection{TV-$L^1$ optical flow: primal-dual and dual form}

By the Fenchel-Rockafellar
duality theory \cite{KK}, under certain regularity condition, the equivalent primal-dual formulation of \eqref{eq:l1:finite} follows 
\begin{equation}\label{eq:primal-dual:optical:o}
\min_{w \in W, \bm{d} \in Y} \max_{p \in Y, q \in Z} \lambda \|\rho(w,\bm{d})\|_{1} + \langle \nabla w, p \rangle 
+ \langle {\nabla \bm d}, q \rangle 
-I_{\{\|p\|_{\infty} \leq 1\}} (p) - I_{\{\|q\|_{\infty} \leq 1\}} (q),
\end{equation}
where $q=(q_1,q_2)^T \in Z= Y \times Y$ and $\nabla \bm d = (\nabla u, \nabla v)^T$.

Actually, the primal-dual form \eqref{eq:primal-dual:optical:o}  is also equivalent to 
\begin{subequations}\label{eq:primal-dual:optical}
	\begin{align}
	\min_{w \in W, u, v \in X} \max_{p \in Y, s \in X, q \in Z}& \langle \nabla w, p \rangle + \langle (\nabla u, \nabla v)^T, q \rangle + \langle I_t + (I_x, I_y) \cdot (u, v)^T + \beta w, s \rangle  \\
	&- I_{\{\|s\|_{\infty} \leq \lambda\}} (s)
	-I_{\{\|p\|_{\infty} \leq 1\}} (p) - I_{\{\|q\|_{\infty} \leq 1\}} (q).
	\end{align}
\end{subequations}
The discrete $L^{\infty}$ and $L^1$ norms are defined as follows \cite{BS, CP}. For $s \in X$, $p = (p_1,p_2)^T \in V$,  $q=(q_1,q_2)^T \in Z$, $q_1 = (q_1^1, q_1^2)^T$, $q_2 = (q_2^1, q_2^2)^T$, $1 \leq t < \infty$,
\begin{subequations}\label{eq:euclid:norm}
	\begin{align}
	&\|s\|_t = \Bigl( \sum_{(i,j) \in \Omega} |s_{i,j}|^t\Bigr)^{1/t}, \quad
	\|s\|_\infty = \max_{(i,j) \in \Omega} \ |s_{i,j}|, \\
	&\|p\|_t = \Bigl( \sum_{(i,j) \in \Omega}
	\bigl( |p_{i,j}|^2 \bigr)^{t/2}\Bigr)^{1/t}, \quad
	\|p\|_\infty = \max_{(i,j) \in \Omega} \ |p_{i,j}|, \quad |p_{i,j}|:= \sqrt{({p_1}_{i,j})^2 + ({p_2}_{i,j})^2}, \\
	& \|q\|_\infty = \max_{(i,j) \in \Omega} \ |q_{i,j}|, \quad |q_{i,j}|: = \sqrt{
		({q_1^1}_{i,j})^2 + ({q_1^2}_{i,j})^2 +  ({q_2^1}_{i,j})^2 + ({q_2^2}_{i,j})^2},\\
	&\|q\|_t = \Bigl( \sum_{(i,j) \in \Omega}
	\bigl( |q_{i,j}|^2)^{t/2}\Bigr)^{1/t}.
	\end{align}
\end{subequations}
Furthermore, by \eqref{eq:dual} and \eqref{eq:dual-cond}, the corresponding dual form of \eqref{eq:l1:finite} is 
\begin{equation}\label{eq:dual:l1:flow:cons}
\max_{p \in Y, s \in X, q \in Z}  \langle I_t, s \rangle 
- I_{\{\|q\|_{\infty} \leq 1\}} (q) - I_{\{\|s\|_{\infty} \leq \lambda\}} (s) -I_{\{\|p\|_{\infty} \leq 1\}} (p),
\end{equation}
subject to the following constraint
\begin{equation}\label{eq:dual:l1:flow:constaint}
\beta s- \Div p=0, \quad \ sI_x - \Div q_1=0, \quad sI_y - \Div q_2=0.
\end{equation}
Similarly, while $\beta=0$, the corresponding primal-dual form of \eqref{eq:l1:finite:no:w} can be written as follows
\begin{subequations}\label{eq:primal-dual:optical:zach}
	\begin{align}
	&\min_{u \in X, v \in X} \max_{ s \in X, q \in Z}  \langle (\nabla u, \nabla v), q \rangle + \langle I_t + (I_x, I_y)^T(u, v), s \rangle - I_{\{\|s\|_{\infty} \leq \lambda\}} (s)  - I_{\{\|q\|_{\infty} \leq 1\}} (q),
	\end{align}
\end{subequations}
and  the corresponding dual form of \eqref{eq:l1:finite:no:w} can be written as
\begin{equation}\label{eq:zach:dual}
\max_{s \in X, q \in Z}  \langle I_t, s \rangle 
- I_{\{\|q\|_{\infty} \leq 1\}} (q) - I_{\{\|s\|_{\infty} \leq \lambda\}} (s),
\end{equation}
subject to the following constraint
\begin{equation}\label{eq:zach:dual:cons}
s[I_x,I_y]'s - \Div q=0, \ \text{i.e.}, \ sI_x - \Div q_1=0, \quad sI_y - \Div q_2=0.
\end{equation}
\subsubsection{Preconditioned ADMM}
With these preparations, we obtain the augmented Lagrangian function for \eqref{eq:dual:l1:flow:cons} and \eqref{eq:dual:l1:flow:constaint},
\begin{subequations}\label{eq:dual:optic}
	\begin{align}
	L_c(w,u,v; p,q,s) =& \langle I_t, s \rangle 
	+ \langle w, \beta s - \D p \rangle +   \langle u,  I_x s - \D q_1 \rangle +  \langle v, I_y s - \D q_2 \rangle \\
	&- \frac{c}{2} \|\beta s - \D p\|^2 - \frac{c}{2} \|I_x s - \D q_1\|^2 - \frac{c}{2} \|I_y s - \D q_2\|^2 \\
	&-I_{\{\|p\|_{\infty} \leq 1\}} (p) -I_{\{\|q\|_{\infty} \leq 1\}} (q) - I_{\{\|s\|_{\infty} \leq \lambda\}} (s). 
	\end{align}
\end{subequations}
The classical augmented Lagrangian method to solve \eqref{eq:dual:optic} is as follows, 
\begin{equation}\label{eq:optic:aug}
\begin{cases}
&(p^{k+1}, q^{k+1}, s^{k+1}) = \text{argmax}_{p,q,s}L_c(w^k,u^k,v^k; p,q,s), \\
&w^{k+1} = w^k - c(\beta s^{k+1} - \D p^{k+1}), \\
&u^{k+1}  = u^k -  c(I_x s^{k+1} - \D q_{1}^{k+1}),\\
&v^{k+1}  = v^k -  c(I_y s^{k+1} - \D q_{2}^{k+1}),
\end{cases}
\end{equation}
where $c$ is the step size.
The main difficulty for applying the augmented Lagrangian method directly as in \eqref{eq:optic:aug} is that it is  very challenging to solve $(p^{k+1}, q^{k+1}, s^{k+1})$ simultaneously due to highly nonlinear and coupling equation of these variables. ADMM is usually considered for solving $(p^{k+1}, q^{k+1}, s^{k+1})$ consecutively. However,
for the subproblem of calculating $(p^{k+1}, q^{k+1}, s^{k+1})$ in  \eqref{eq:optic:aug} by ADMM, there is no convergence guarantee by solving it consecutively since there are 3 block of variables \cite{CBHY}. Now, we will circumvent this difficulty by regrouping the 3 block of variables into 2 block of variables.  Actually, with notation $y:=(p,q)^T$,   ${\Lambda}: = (w,u,v)$, and 
\[
\mathcal{A} = (\beta I, I_x, I_y)^T, \quad \mathcal{B}=\text{Diag}[-\D, -\D], \quad \mathcal{B}^* = \text{Diag}[\nabla, \nabla],
\]
we can rewrite equation \eqref{eq:dual:optic} as follows 
\begin{equation}\label{eq:dual:optic:two:block}
L_c(\Lambda; y,s) = \langle I_t, s \rangle +  \langle \Lambda, \mathcal{A}s+\mathcal{B}y \rangle
- \frac{c}{2} \|\mathcal{A}s+\mathcal{B}y \|^2 
-\mathcal{ G} (s) -\mathcal{H}(y),
\end{equation}
where  $\mathcal{G}(s) =  I_{\{\|s\|_{\infty} \leq 1\}} (s) - \langle I_t, s \rangle  $ and $\mathcal{H}(y) = I_{\{\|q\|_{\infty} \leq 1\}} (q) + I_{\{\|p\|_{\infty} \leq 1\}} (p)$.
Therefore, we can use the preconditioned or semi-proximal ADMM as follows
\begin{subequations}\label{eq:aug:pre1:vector}
	\begin{align}
	y^{k+1} &=  \text{argmax}_{y}L_c(\Lambda^k; y,s^k) - \frac{1}{2}\|y-y^k\|_{acI-c\nabla \nabla^*}^2,\\
	s^{k+1} &=  \text{argmax}_{s}L_c(\Lambda^k; y^{k+1},s) - \frac{1}{2}\|s-s^k\|_{\tilde acI-cM_s}^2,\\
	\Lambda^{k+1} &= \Lambda^k-rc(\mathcal{A}s^{k+1}+\mathcal{B}y^{k+1} ),
	\end{align}
\end{subequations}
where $M_s = \beta^2 + I_x^2 + I_y^2$ and $r \in (0, \frac{\sqrt{5}+1}{2})$ is the relaxation parameter \cite{DY, LST}. This kind of relaxation that is only  relaxed on the update of the Lagrangian multipliers is originated from \cite{FG}. The weighted $L^2$ norm was introduced in \cite{DY, LST} where the corresponding weight should be positive semidefinite for the convergence. Thus the specially chosen $a$ or $\tilde a$ can be as follows,
\begin{align}
&\|z-z^k\|_{M}^2 = \langle (z-z^k),M(z-z^k)\rangle, \quad z = y \  \text{or} \ s, \\
& acI-c \nabla \nabla^*  \geq 0, \quad \tilde acI-c M_s \geq 0. 
\end{align}
We thus choose $a$ and $\tilde a$ as follows
\begin{equation}\label{eq:stepsize:constraint}
\frac{1}{a} \leq \frac{1}{\|\nabla \nabla^*\|}, \quad \frac{1}{\tilde a}  \leq \frac{1}{\|M_s\|}.
\end{equation}
Writing \eqref{eq:aug:pre1:vector} component-wisely, we have the relaxed and preconditioned ADMM
\begin{equation}\label{eq:rpADMMI:cp}\tag{rpADMMI}
\begin{cases}
p^{k+1} =  \mathcal{P}_p[(I - \frac{1}{a}\nabla \nabla^*)p^k +\frac{1}{a}(\frac{\nabla w^k}{c} - \beta \nabla s^k)] ,\\
q_{1,t}^{k+1} = [(I - \frac{1}{a}\nabla \nabla^*)q_1^k +\frac{1}{a}(\frac{\nabla u^k}{c} -  \nabla (s^kI_x))], \\
q_{2,t}^{k+1} = [(I - \frac{1}{a}\nabla \nabla^*)q_2^k +\frac{1}{a}(\frac{\nabla v^k}{c} -  \nabla (s^kI_y))], \\
(q_1^{k+1}, q_2^{k+1}) = \mathcal{P}_q[(q_{1,t}^{k+1}, q_{2,t}^{k+1})],\\
s^{k+1} = \mathcal{P}_s[(I- \frac{1}{\tilde a}M_s)s^k + \frac{1}{\tilde ac}(I_t + \beta w^k + u^k I_x + v^k I_y \\
\quad \quad \quad \ + c\beta \D p^{k+1} + c I_x \D q_1^{k+1} + c I_y \D q_2^{k+1}) ], \\
w^{k+1} = w^k - rc(\beta s^{k+1} - \D p^{k+1}), \\
u^{k+1}  = u^k -  rc(I_x s^{k+1} - \D q_{1}^{k+1}),\\
v^{k+1}  = v^k -  rc(I_y s^{k+1} - \D q_{2}^{k+1}),
\end{cases}
\end{equation}
where the projections are defined as follows \cite{CP, BS}
\begin{align*}
&\mP_s(s) = (I + \frac{1}{\tilde a c} \partial \mathcal{G})^{-1}(s)
= \argmin_{s' \in X} \ \frac12 \| s' - s\|^2_2
+ \frac{1}{\tilde a c} \mathcal{I}_{\{\|s\|_\infty \leq 1\}}(s') 
= \frac{s}{\max(1, |s|/\alpha)}, \\
&\mP_p(p) = (I + \frac{ \partial I_{\{\|p\|_{\infty} \leq 1\}}}{ a c})^{-1}(p)
= \argmin_{p' \in X} \ \frac12 \| p' - p\|^2_2
+ \frac{\mathcal{I}_{\{\|p\|_\infty \leq 1\}}(p') }{\tilde a c} 
= \frac{p}{\max(1, |p|/\alpha)}, \\
&\mP_q(q) = (I + \frac{ \partial I_{\{\|q\|_{\infty} \leq 1\}}}{ a c})^{-1}(q)
= \argmin_{q' \in X} \ \frac12 \| q' - q\|^2_2
+ \frac{\mathcal{I}_{\{\|q\|_\infty \leq 1\}}(q') }{\tilde a c} 
= \frac{q}{\max(1, |q|/\alpha)},
\end{align*}
where $|s|$, $|p|$ and $|q|$ are defined in \eqref{eq:euclid:norm}. For the detail of the calculations, we refer to the Appendix \ref{sec:appendix}.

There is another kind of relaxation of preconditioned ADMM. Unlike the relaxation in \eqref{eq:aug:pre1:vector}, the relaxed and preconditioned ADMM for solving the dual problem \eqref{eq:dual:optic} with the augmented Lagrangian \eqref{eq:dual:optic:two:block} reads as follows \cite{SUN},
\begin{align}
y^{k+1} &=(N + \partial \mathcal{H})^{-1}( \mathcal{B}^*(-c \mathcal{A}s^k + \Lambda^k) + (N-c \mathcal{B}^* \mathcal{B})y^k), \notag \\
s^{k+1} &= (M + \partial  \mathcal{G})^{-1}( \mathcal{A}^*(-c\rho  \mathcal{B}y^{k+1} +c(1-\rho ) \mathcal{A}s^k + \Lambda^k)+(M-r \mathcal{A}^* \mathcal{A})s^k), \label{eq:update:two:sun:lag} \\
\Lambda^{k+1} &=\Lambda^k - r( \mathcal{A}s^{k+1}-(1-\rho ) \mathcal{A}s^k+\rho  \mathcal{B}y^{k+1}). \notag
\end{align}
Here $\rho \in (0,2)$ is the relaxation parameter  \cite{EP,SUN} and $N$, $M$ are two linear and bounded operators \cite{SUN}, such that
\[
N-c \mathcal{B}^* \mathcal{B} \geq 0 \Rightarrow N - c \nabla \nabla^* \geq 0 , \quad M-c \mathcal{A}^* \mathcal{A} \geq 0.
\]
We thus choose $N = acI$, $M=\tilde a cI$ with $a$, $\tilde a$ satisfy the condition \eqref{eq:stepsize:constraint}.
This kind of relaxation is originated from the relaxed Douglas-Rachford splitting method \cite{EP}, since one can get the relaxed ADMM by applying the Douglas-Rachford splitting method to the dual problem \cite{EP,SUN}. Preconditioned techiques are introduced for both blocks with mild conditions \cite{BS0,SUN}. Writing \eqref{eq:update:two:sun:lag}   component-wisely with application to \eqref{eq:dual:optic:two:block}, we have 
\begin{equation}\label{eq:rpADMMII:sun}\tag{rpADMMII}
\begin{cases}
p^{k+1} =  \mathcal{P}_p[(I - \frac{1}{a}\nabla \nabla^*)p^k +\frac{1}{a}(\frac{\nabla w^k}{c} - \beta \nabla s^k)] ,\\
q_{1,t}^{k+1} = [(I - \frac{1}{a}\nabla \nabla^*)q_1^k +\frac{1}{a}(\frac{\nabla u^k}{c} -  \nabla (s^kI_x))], \\
q_{2,t}^{k+1} = [(I - \frac{1}{a}\nabla \nabla^*)q_2^k +\frac{1}{a}(\frac{\nabla v^k}{c} -  \nabla (s^kI_y))], \\
(q_1^{k+1}, q_2^{k+1}) = \mathcal{P}_q{(q_{1,t}^{k+1}, q_{2,t}^{k+1})},\\
s^{k+1} = \mathcal{P}_s\big\{(I- \frac{\rho}{\tilde a}M_s)s^k + \frac{1}{\tilde ac}[I_t + \beta w^k + u^k I_x + v^k I_y \\
\quad \quad \quad \ + c \rho(\beta \D p^{k+1} +  I_x \D q_1^{k+1} + I_y \D q_2^{k+1})] \big\}, \\
w^{k+1} = w^k - c(\beta s^{k+1} - \rho \D p^{k+1} - (1-\rho)\beta s^{k}), \\
u^{k+1}  = u^k -  c(I_x s^{k+1} - \rho\D q_{1}^{k+1} - (1-\rho)I_x s^{k}),\\
v^{k+1}  = v^k -  c(I_y s^{k+1} - \rho \D q_{2}^{k+1}- (1-\rho)I_y s^{k}).
\end{cases}
\end{equation}
Similarly, the augmented Lagrangian function for \eqref{eq:zach:dual} and \eqref{eq:zach:dual:cons} becomes
\begin{subequations}\label{eq:dual:optic:zach}
	\begin{align}
	L_c(u,v; q,s) =& \langle I_t, s \rangle 
	+   \langle u,  I_x s - \D q_1 \rangle +  \langle v, I_y s - \D q_2 \rangle - I_{\{\|q\|_{\infty} \leq 1\}} (q) \\
	& - I_{\{\|s\|_{\infty} \leq \lambda\}} (s)- \frac{c}{2} \|I_x s - \D q_1\|^2 - \frac{c}{2} \|I_y s - \D q_2\|^2.
	\end{align}
\end{subequations}
Since there are only two block of variables, we can directly use the preconditioned ADMM
\begin{subequations}\label{eq:aug:pre1:vector:zach}
	\begin{align}
	q^{k+1} &=  \text{argmax}_{y}L_c(u^k,v^k; q,s^k) - \frac{1}{2}\|q-q^k\|_{acI-c\nabla^*\nabla}^2,\\
	s^{k+1} &=  \text{argmax}_{s}L_c(u^k,v^k; q^{k+1},s) - \frac{1}{2}\|s-s^k\|_{\tilde acI-cM_s}^2,
	\end{align}
\end{subequations}
where the constraints on $a$ and $\tilde a$ are still chosen by   \eqref{eq:stepsize:constraint} except $\beta=0$ here. Writing \eqref{eq:aug:pre1:vector:zach} component-wisely and together with the updates of the Lagrangian multipliers, i.e., $u^{k+1}$ and $v^{k+1}$, we arrive at preconditioned ADMM for the original TV-$L^1$ optical flow estimate \eqref{eq:l1:finite:no:w} without relaxation
\begin{equation}\label{eq:admm:zach}\tag{Zach-pADMM}
\begin{cases}
q_{1,t}^{k+1} = [(I - \frac{1}{a}\nabla \nabla^*)q_1^k +\frac{1}{a}(\frac{\nabla u^k}{c} -  \nabla (s^kI_x))], \\
q_{2,t}^{k+1} = [(I - \frac{1}{a}\nabla \nabla^*)q_2^k +\frac{1}{a}(\frac{\nabla v^k}{c} -  \nabla (s^kI_y))], \\
(q_1^{k+1}, q_2^{k+1}) = \mathcal{P}_q[(q_{1,t}^{k+1}, q_{2,t}^{k+1})],\\
s^{k+1} = \mathcal{P}_s[(I- \frac{1}{\tilde a}M_s)s^k + \frac{1}{\tilde ac}(I_t + u^k I_x + v^k I_y \\
\quad \quad \quad \ + c\beta \D p^{k+1} + c I_x \D q_1^{k+1} + c I_y \D q_2^{k+1}) ], \\
u^{k+1}  = u^k -  c(I_x s^{k+1} - \D q_{1}^{k+1}),\\
v^{k+1}  = v^k -  c(I_y s^{k+1} - \D q_{2}^{k+1}).
\end{cases}
\end{equation}
For the convergence of \ref{eq:rpADMMI:cp}, \ref{eq:rpADMMII:sun} and \ref{eq:admm:zach}, we have the following proposition \cite{DY, LST, BS0, SUN}.
\begin{proposition}\label{thm:twolabel}
	For the relaxed and preconditioned ADMM of type \ref{eq:rpADMMI:cp},  if choosing $a=8$ and $\tilde a = \|M_s\|$, we can get the convergence  of the iteration \ref{eq:rpADMMI:cp} for any $r \in (0, \frac{\sqrt{5}+1}{2})$.  For the relaxed and preconditioned ADMM of type \ref{eq:rpADMMII:sun}, if  choosing  $M=c \|M_s\|I$ and $N=8cI$ as in \eqref{eq:update:two:sun:lag}, we get the convergence of the iteration \ref{eq:rpADMMII:sun} for any $\rho \in (0,2)$. For both \ref{eq:rpADMMI:cp} and \ref{eq:rpADMMII:sun}, the dual sequence $( p^k, q^k, s^k)$ converges to the solution $( p^*, q^*, s^*)$ of the dual problem \eqref{eq:dual:l1:flow:cons} and the Lagrangian multipliers $( w^k, u^k, v^k)$ converge to the solution $( w^*, u^*, v^*)$  of the primal problem \eqref{eq:l1:finite}. The corresponding  ergodic convergence rate of the primal and dual iteration sequences of \ref{eq:rpADMMI:cp} or \ref{eq:rpADMMII:sun} is $\mathcal{O}(1/k)$.
\end{proposition}

\begin{remark}
	The convergence of \ref{eq:admm:zach} for solving \eqref{eq:l1:finite:no:w} is completely similar to Proposition \ref{thm:twolabel} including the convergence of iteration sequence and the corresponding ergodic convergence rate \cite{DY,LST,SUN}.
\end{remark}
\subsubsection{Preconditioned Douglas-Rachford splitting method}
Now, let's turn to  the preconditioned Douglas-Rachford splitting algorithm for \eqref{eq:primal-dual:optical:o}. The preconditioned Douglas-Rachford splitting algorithm is an efficient algorithm aiming at dealing with  the  challenging linear subproblems during each nonlinear Douglas-Rachford iteration with any finite feasible preconditioned iterations  while applying the method for solving nonlinear saddle-point problems; see \cite{BS} for its development and applications in image restoration problems.  With the data
$F(\mathfrak{X}) = \lambda \|\rho(w,\bm d) \|_{1}$, $G(y) = I_{\{\|p\|_{\infty} \leq 1\}} (p) + I_{\{\|q\|_{\infty} \leq 1\}} (q)$ and $\mathcal{K} = \text{Diag}(\nabla, \nabla)$ and notation $\mathfrak{X}:=(w, \bm d)^T$, we can reformulate \eqref{eq:primal-dual:optical:o} as the following generic  saddle-point problem 
\begin{equation}\label{eq:saddle:qr}
\min_{\mathfrak{X}} \max_{y} {F}(\mathfrak{X}) + \langle \mathcal{K}\mathfrak{X},y \rangle -{G}(y),
\end{equation}
where $y=(p,q)^T$ as before.

Each iteration of the  preconditioned  Douglas-Rachford splitting method  for \eqref{eq:saddle:qr} can be written as:
\begin{equation}\label{iteration:general:relax:dual} \tag{pDR}
\begin{cases}
b^{k}  = \bar {\mathfrak{X}}^k -\sigma K^* \bar y^k, \\
\mathfrak{X}^{k+1} = \mathfrak{X}^k+ M^{-1}(b^{k} - T \mathfrak{X}^k), \\
y^{k+1} = \bar y^{k} + \tau \mathcal{K} \mathfrak{X}^{k+1}, \\
\bar {\mathfrak{X}}^{k+1} = \bar {\mathfrak{X}}^k + [(I + \sigma \partial F)^{-1}
[2  \mathfrak{X}^{k+1} - \bar {\mathfrak{X}}^k] -  \mathfrak{X}^{k+1}], \\
\bar{y}^{k+1} = \bar y^k + [( I + \tau \partial G)^{-1}
[2y^{k+1} - \bar y^k] - y^{k+1}],
\end{cases}
\end{equation}
where  $\sigma$, $\tau$ are positive step sizes that can be chosen freely \cite{BS, BS3}. $M$ is the feasible preconditioner for $T = I + \sigma \tau \mathcal{K}^*\mathcal{K}$ \cite{BS}, i.e., 
\[
M - T \geq 0 \Leftrightarrow M-T \  \text{is positive semidefinite.}
\]
The convergence of iterations \eqref{iteration:general:relax:dual} can be guaranteed \cite{BS,BS3} with the above feasibility condition. Supposing $ \mathfrak{X} = ( w, {\bm d})$ and  $\tilde {\mathfrak{X}} = (\tilde w, \tilde {\bm d})$, 
the resolvent $\mathfrak{X}=(I + \sigma \partial F)^{-1}(\tilde {\mathfrak{X}})$ can be found in \cite{CP} or \cite{ZPH}, 
\begin{equation}\label{eq:resolvent}
(w_{i,j}, \bm d_{i,j}) = (\tilde w_{i,j},\tilde {\bm d}_{i,j}) +
\begin{cases}
\sigma \lambda e_{i,j}& \quad \rho(\tilde w_{i,j},\tilde x_{i,j}) < -\sigma \lambda |e|_{i,j}^2 \\
-\sigma \lambda e_{i,j}& \quad \rho(\tilde w_{i,j},\tilde x_{i,j}) > \sigma \lambda |e|_{i,j}^2 \\
-\rho(\tilde w_{i,j},\tilde x_{i,j})e_{i,j}/|e|_{i,j}^2& \quad |\rho(\tilde w_{i,j},\tilde x_{i,j})| \leq \sigma \lambda |e|_{i,j}^2 
\end{cases}
\end{equation}
where $e_{i,j} = (\beta, (\nabla I)_{i,j})$ and $|e|_{i,j}^2 = \beta^2 + |\nabla I|_{i,j}^2$.
For the convenience of the reader, we give an elementary proof with the Fenchel-Rockafellar duality in the Appendix \ref{sec:appendix}.

Now, we will employ the preconditioned Douglas-Rachford splitting framework as in \cite{BS} for \eqref{eq:saddle:qr}, where the classical prconditioned iterations for linear equation is cooperated into the nonlinear Douglas-Rachford splitting method and the convergence can be guaranteed without error control as in \cite{EP}.
The preconditioned Douglas-Rachford splitting algorithm for \eqref{eq:primal-dual:optical:o} with saddle-point structure \eqref{eq:saddle:qr} is as follows,
\begin{equation}\label{iteration:opflow:l1:pdr}\tag{pDR}
\begin{cases}
w^{k+1} = w^{k} + M^{-1}[ \bar w^k + \sigma \D p^k- T w^{k}], \\
u^{k+1} = u^{k} + M^{-1}[ \bar u^k + \sigma \D q_{1}^k- T u^{k}], \\
v^{k+1} = v^{k} + M^{-1}[ \bar v^k + \sigma \D q_{2}^k- T v^{k}], \\
(\bar w^{k+1},\bar x^{k+1})^{T} = (\bar w^{k},\bar x^{k})^{T}
+  [(I + \sigma \partial {F})^{-1}
[2 ( w^{k+1}, x^{k+1})^{T} \\
\qquad \quad \quad \quad \quad \quad \quad - (\bar w^{k},\bar x^{k})^{T}] - ( w^{k+1}, x^{k+1})^{T}], \\
\bar p^{k+1} = -  \tau \nabla w^{k+1} + [\mathcal{P}_{p}(\bar p^k + 2 \tau \nabla w^{k+1} )],\\
\bar q^{k+1} = -  \tau (\nabla u^{k+1}, \nabla v^{k+1})^{T} + [\mathcal{P}_{q}(\bar q^k + 2 \tau (\nabla u^{k+1}, \nabla v^{k+1})^{T} )],
\end{cases}
\end{equation}
where $M$ is any finite times classical symmetric Gauss-Seidel preconditioner which is a feasible preconditioner for $ T:=I -\sigma \tau \Delta$ with Neumann boundary condition \cite{BS}. The symmetric Gauss-Seidel iterations are designed for dealing with the following perturbed Laplacian equation of $u$ ($w$ or $v$)      
\begin{equation}\label{neumann:boudary}
\left\{
\begin{aligned}
u - \sigma \tau \Delta u &= b^k,\\
\frac{\partial u}{\partial \nu}|_{\partial \Omega} &=  0,
\end{aligned}
\right.
\end{equation}%
which appeared in the original Douglas-Rachford splitting method.
Actually, any finite symmetric Red-Black Gauss--Seidel method (SRBGS)  can be formulated as the following classical preconditioned iteration
\[
u^{k+1}: = u^k + M^{-1}(b^k- T u^k), \ \text{with} \ M \geq T \Leftrightarrow   \text{M is a feasible preconditioner for} \ T.
\]
With SRBGS, we can get the more efficient and convergent preconditioned Douglas-Rachford splitting iterations as in \ref{iteration:opflow:l1:pdr}. 
For the convergence of the preconditioned Douglas-Rachford splitting method \ref{iteration:opflow:l1:pdr}, we have the following proposition \cite{BS,BS3}.

\begin{proposition}[\cite{BS}]\label{pro:pdr}
	If the preconditioner satisfies the feasibility condition, i.e.,  $M \geq T$, then iteration sequence $\{   w^k, u^k, v^k, p^k, q^k, s^k\}$ of   the preconditioned Douglas-Rachford splitting \ref{iteration:opflow:l1:pdr}  converges to a saddle-point $(  w^*, u^*, v^*, p^*, q^*, s^*)$ of \eqref{eq:primal-dual:optical:o} with $( w^*, u^*, v^*)$ being a solution of the primal problem \eqref{eq:l1:finite} and $(p^*, q^*, s^*)$ being a solution of the dual problem \eqref{eq:dual:l1:flow:cons}.	The ergodic convergence rate of the iteration sequences $\{   w^k, u^k, v^k, p^k, q^k, s^k\}$ is $\mathcal{O}(1/k)$.
\end{proposition}

\section{Numerical experiments}\label{sec:num}
In this part, we will study the numerical performance of the proposed preconditioned ADMM algorithms. For the optical flow estimates, we integrate the algorithms into a standard  coarse-to-fine framework, which is an efficient multiscale pyramid process \cite{PBS} that can greatly reduce the optimization energy and improve the performance. 
The algorithms for comparison and the corresponding parameter settings  are as follows:
\begin{itemize}
	\item For ALG1 with application to \eqref{eq:l1:finite},  
	the primal-dual algorithm
	introduced in \cite{CP} with constant step sizes: the dual step size $\tau=0.1$,  $\sigma = 1/(L^2\sigma)$ with $L = \sqrt{8}$.
	\item For pADMM, \ref{eq:rpADMMI:cp} and \ref{eq:rpADMMII:sun} with application to \eqref{eq:l1:finite} : we choose  $c=0.05$. pADMM is the \ref{eq:rpADMMI:cp} without relaxation, i.e., $r=1.0$. We choose the relaxation parameters $r=1.618$ for \ref{eq:rpADMMI:cp} and $\rho=1.9$ for \ref{eq:rpADMMII:sun}.
	\item For \ref{iteration:opflow:l1:pdr} with application to \eqref{eq:l1:finite}: we choose $\sigma = 2$, $\tau = 0.4$. 
	\item For \ref{eq:admm:zach} with application to \eqref{eq:l1:finite:no:w}: we choose $c=0.05$.
	\item Zach-O: we use it to denote the alternating minimization method for \eqref{eq:alter:appro} with the default parameter settings in \cite{PMF} through the corresponding online demo  of Image Processing On Line with manually loaded data \cite{PMF}. Note that the anisotropic model \eqref{eq:alter:appro} is different from the isotropic model \eqref{eq:l1:optical:flow:no:w}.  
	
\end{itemize}
We use the Middlebury optical flow benchmark data set with ground truth (\url{http://vision.middlebury.edu/flow/}). We also employ the Ettlinger--Tor and the Rheinhafen sequences (\url{http://www.ira.uka.de/image_sequences/}). The Yosemite test sequences are also used for comparison. For \eqref{eq:l1:finite}, we choose $\beta=0.05$ or $\beta=0.01$.

Table \ref{tab:optical:zachp:beta} first shows the comparison between the model \eqref{eq:l1:finite} and \eqref{eq:l1:finite:no:w} with different algorithms. For Table \ref{tab:optical:zachp:beta}, we employ the same step size $c$ and multiscale pyramid setting for pADMM, \ref{eq:rpADMMI:cp}, \ref{eq:rpADMMII:sun} and \ref{eq:admm:zach}.  It can be seen the model \eqref{eq:l1:finite} performs better than the model \eqref{eq:l1:finite:no:w} with lower average angular error and end point error. pADMM and pDR can give the lowest average angular error. The over-relaxation \ref{eq:rpADMMI:cp} and \ref{eq:rpADMMII:sun} do not bring out better results except the Hydreangea with \ref{eq:rpADMMII:sun}. Actually, it was  proved theoretically that the ergodic convergence rate can be better with over-relaxation \cite{DY, LST, SUN} and it was also shown numerically one can get faster algorithm with over-relaxation \cite{SUN, SYT}. However, there is no promising improvement here for optical flow estimates, compared pADMM with \ref{eq:rpADMMI:cp} or \ref{eq:rpADMMII:sun} for the same model \eqref{eq:l1:finite}. It is probably because of the multiscale pyramid structure, which is nonlinear and the initial values are continuously changed for different scales. 

\begin{figure}
	\begin{center}
		\subfloat[RubberWhale: frame10]
		{\includegraphics[width=0.28\textwidth]{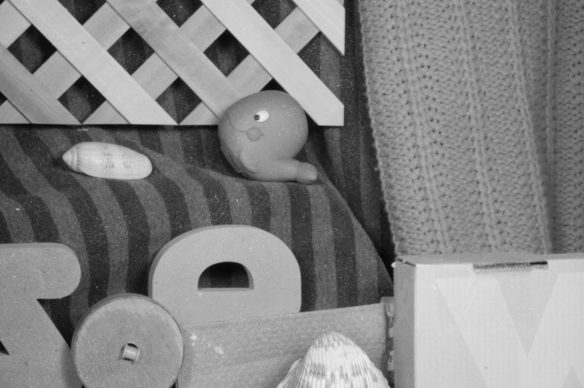}}\ \ \ 
		\subfloat[RubberWhale: flow, pADMM]
		{\includegraphics[width=0.28\textwidth]{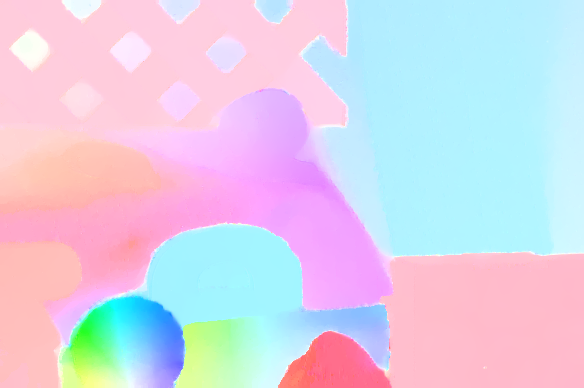}} \ \ \ 
		\subfloat[RubberWhale: flow, Zach-O]
		{\includegraphics[width=0.28\textwidth]{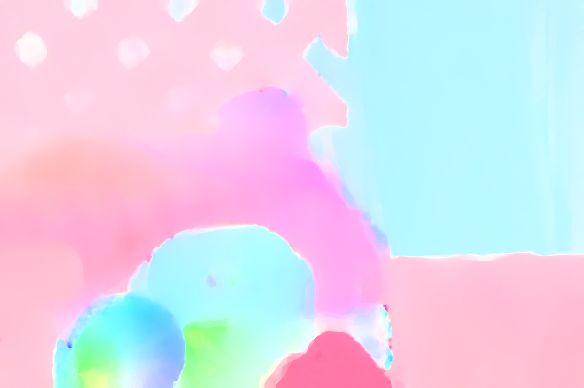}}	
		\\ [-0.1em]
		\subfloat[Urban2: frame10]
		{\includegraphics[width=0.28\textwidth]{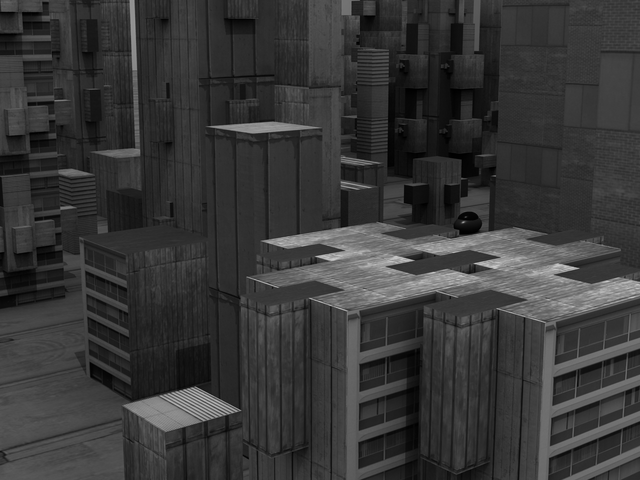}}\ \ \ 
		\subfloat[Urban2: flow, pADMM]
		{\includegraphics[width=0.28\textwidth]{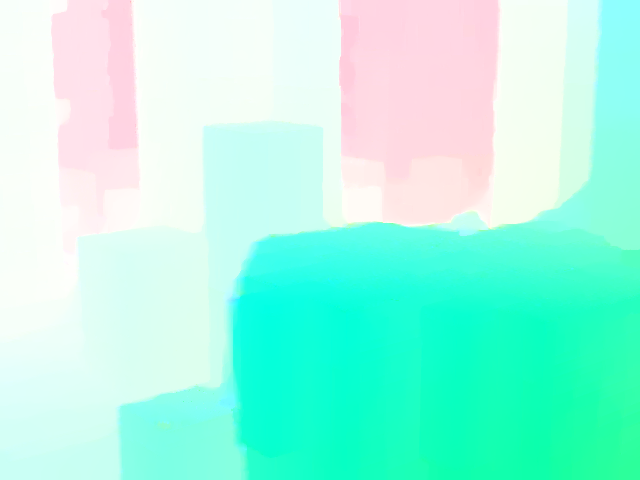}}\ \ \ 
		\subfloat[Urban2: flow, Zach-O]
		{\includegraphics[width=0.28\textwidth]{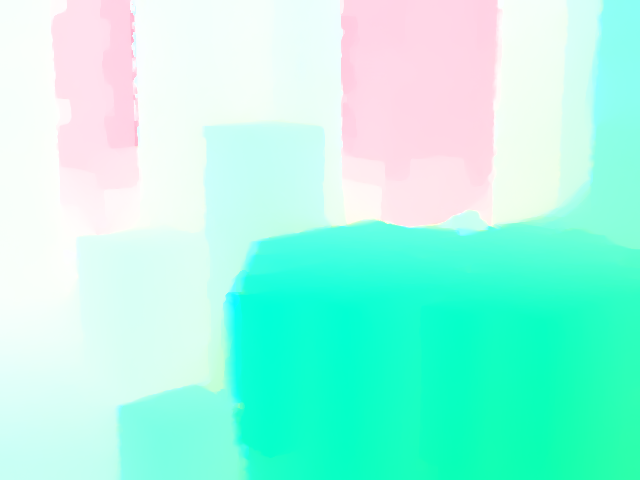}}\\ [-0.1em]
		\subfloat[Grove2: frame10]
		{\includegraphics[width=0.28\textwidth]{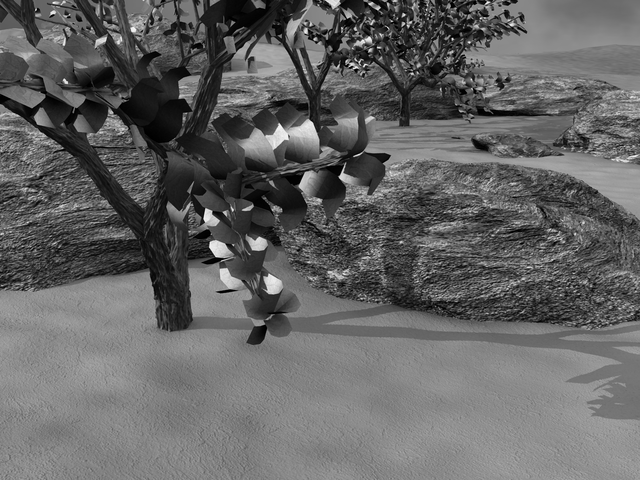}}\ \ \ 
		\subfloat[Grove2: flow, pADMM]
		{\includegraphics[width=0.28\textwidth]{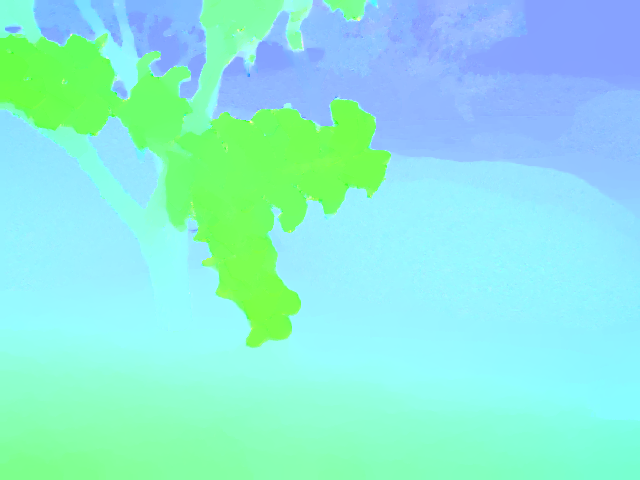}}\  \ \ 
		\subfloat[Grove2: flow, Zach-O]
		{\includegraphics[width=0.28\textwidth]{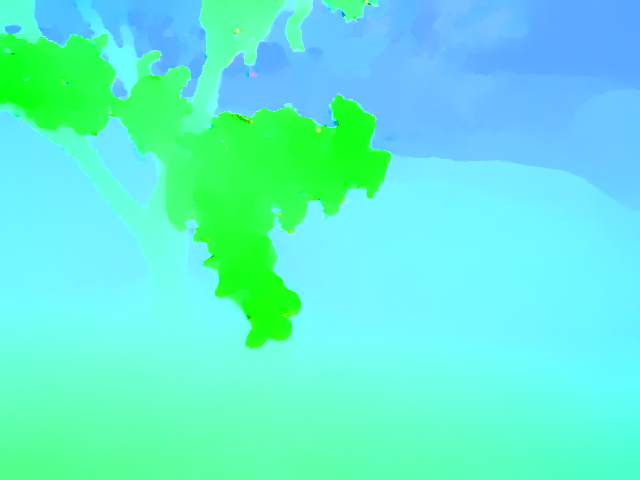}}\\ [-0.1em]
		\subfloat[Hydrangea: frame10]
		{\includegraphics[width=0.28\textwidth]{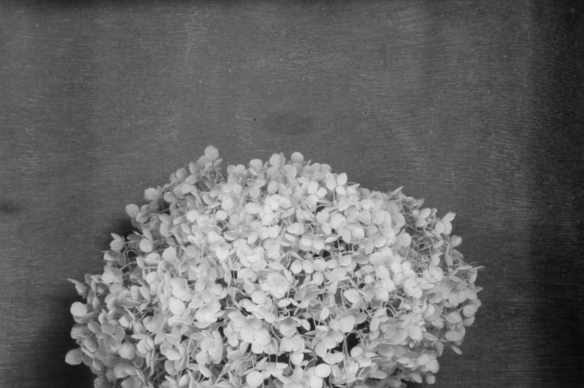}}\ \ \ 
		\subfloat[Hydrangea: flow, pADMM]
		{\includegraphics[width=0.28\textwidth]{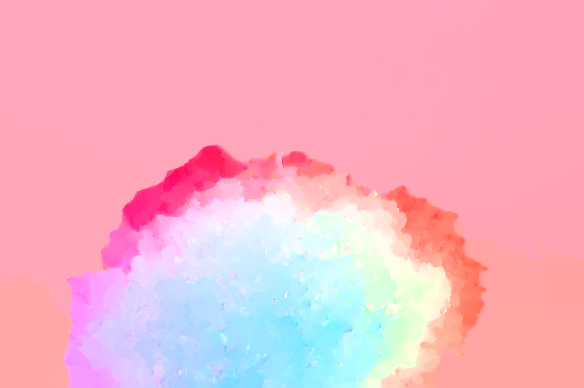}}\ \ \ 
		\subfloat[Hydrangea: flow, Zach-O]
		{\includegraphics[width=0.28\textwidth]{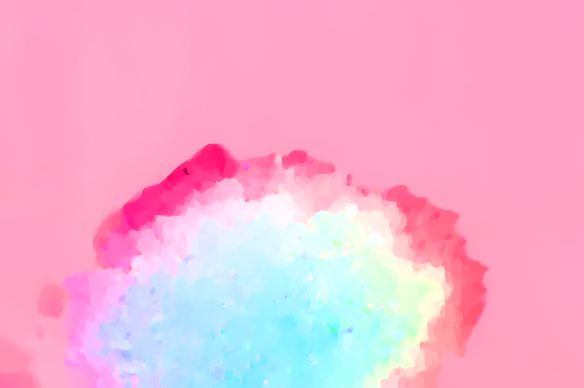}}\\ [-0.1em]
		\subfloat[Venus: frame10]
		{\includegraphics[width=0.28\textwidth]{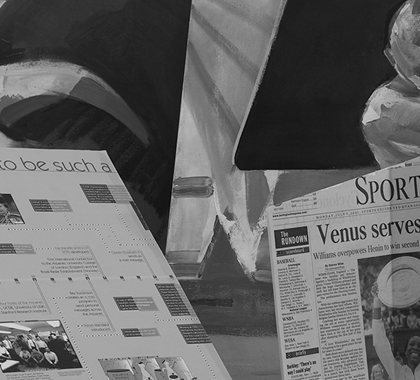}}\ \ \ 
		\subfloat[Venus: flow, pADMM]
		{\includegraphics[width=0.28\textwidth]{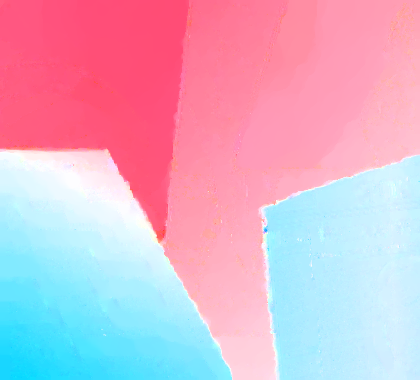}}\ \ \ 
		\subfloat[Venus: flow, Zach-O]
		{\includegraphics[width=0.28\textwidth]{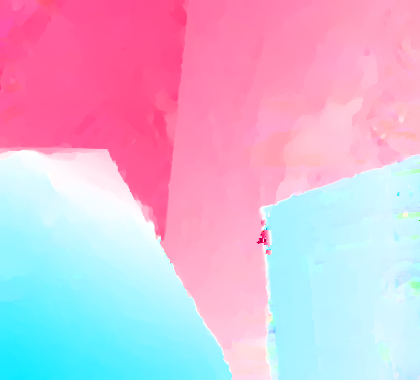}}
	\end{center}
	\vspace*{-0.5em}
	\caption{The results of the optical flow estimates. The images on the left column are the Frame 10 image from the Middlebury test sequences. The images in the middle column are the optical flow estimates by pADMM. The images on the right column are the optical flow estimates by Zach-O.}
	\label{optical:with:gt}
\end{figure}

Figure \ref{optical:with:gt} shows the reconstructed optical flow estimates with the proposed preconditioned ADMM algorithms. pADMM can provide more sharp and clean reconstructed  optical flow estimate than Zach-O for the Venus and RubberWhale sequences. This is probably because the model \eqref{eq:l1:optical:flow:zach} is anisotropic on $\bm d$ while the model \eqref{eq:l1:optical:flow:no:w} is isotropic on $\bm d$ and  Zach-O algorithm is essentially based on the approximate model of \eqref{eq:alter:appro} for the original model \eqref{eq:l1:optical:flow:zach}.

Table \ref{tab:angular:error1} shows the comparison between the Zach-O algorithm and the proposed preconditioned ADMM algorithms including the pADMM, \ref{eq:rpADMMI:cp}, \ref{eq:rpADMMII:sun}. pADMM, \ref{eq:rpADMMI:cp} or \ref{eq:rpADMMII:sun} give better 
average angular error and average end point error compared with the results with best parameters by Zach-O \cite{PMF} (see Table 4 in \cite{PMF}) except the Venus sequence.

\begin{figure}
	\begin{center}   
		\subfloat[Yosemite: frame 08]
		{\includegraphics[width=0.30\textwidth]{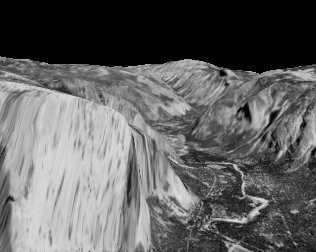}}\ \ \ 
		\subfloat[Yosemite: frame 09]
		{\includegraphics[width=0.30\textwidth]{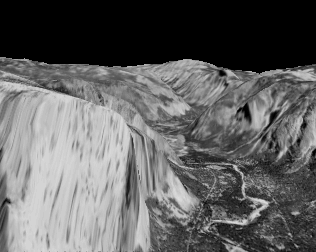}}\ \ \ 
		\subfloat[Yosemite: flow]
		{\includegraphics[width=0.30\textwidth]{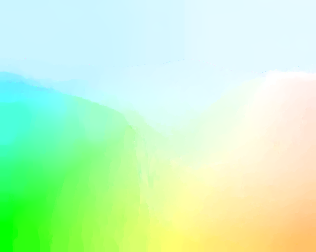}}\\ [-0.1em]	
		\subfloat[Ettlinger--Tor: frame 07]
		{\includegraphics[width=0.30\textwidth]{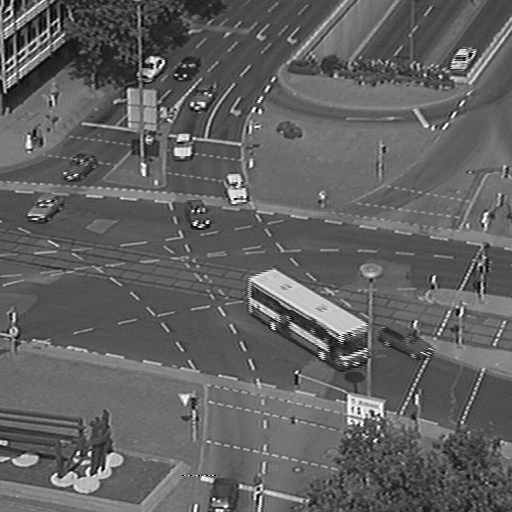}}\ \ \ 
		\subfloat[Ettlinger--Tor: frame 08]
		{\includegraphics[width=0.30\textwidth]{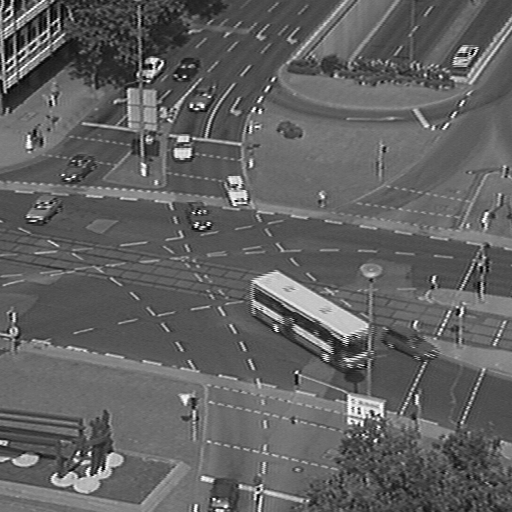}}\ \ \ 
		\subfloat[Ettlinger Tor: flow]
		{\includegraphics[width=0.30\textwidth]{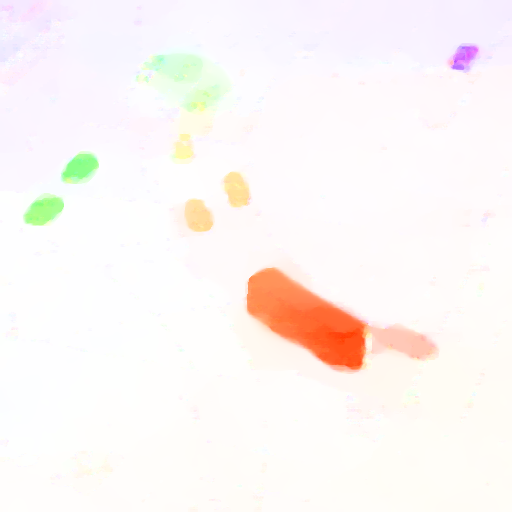}}\\ [-0.1em]
		\subfloat[Rheinhafen: frame 1130]
		{\includegraphics[width=0.30\textwidth]{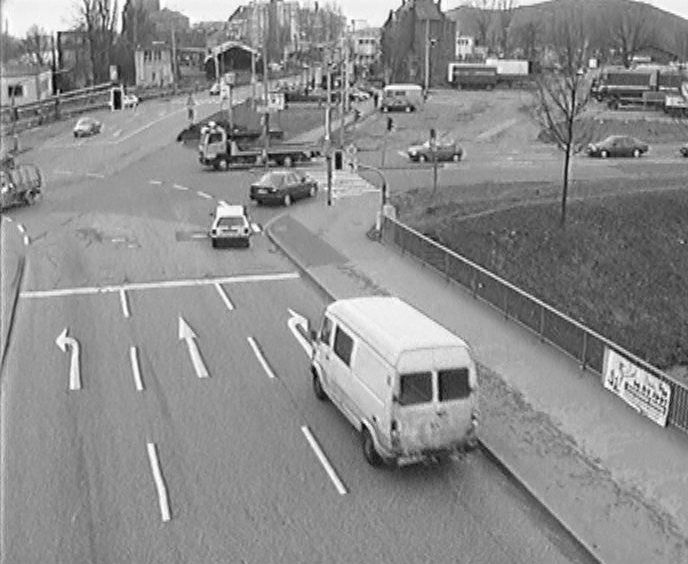}}\ \ \ 
		\subfloat[Rheinhafen: frame 1131]
		{\includegraphics[width=0.30\textwidth]{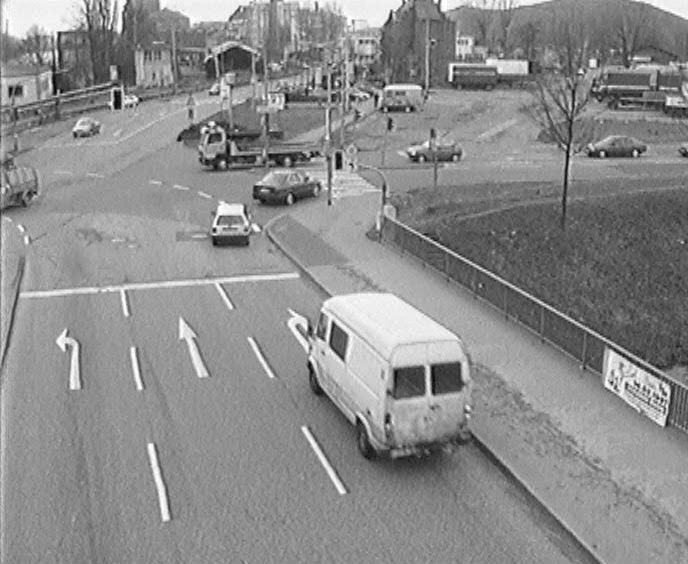}}\ \ \ 
		\subfloat[Rheinhafen: flow]
		{\includegraphics[width=0.30\textwidth]{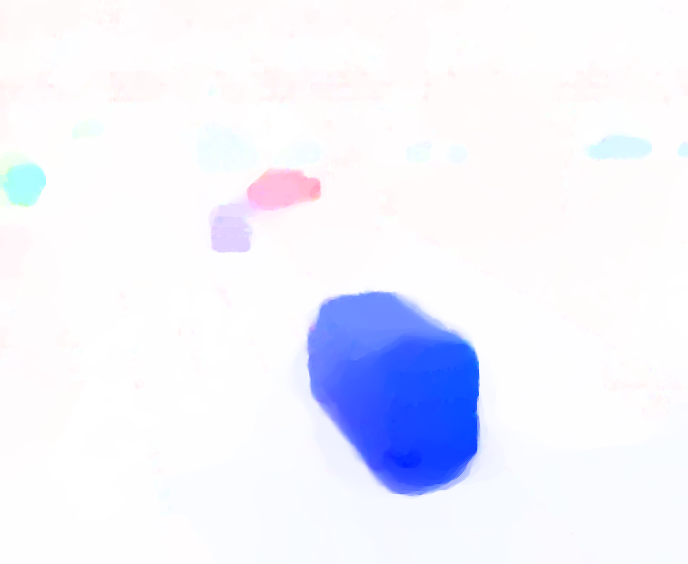}}
	\end{center}
	\vspace*{-0.5em}
	\caption{The results of the optical flow estimates. The images on the left and the middle  columns are the corresponding image sequences. The images on the right column are the optical flow estimates by pADMM.}
	\label{yosemite:erheinhafen}
\end{figure}

Figure \ref{yosemite:erheinhafen} shows the optical flow estimates on Yosemite, Ettlinger--Tor and Rheinhafen sequences.  The corresponding optical flow estimates of pADMM are also of high quality.

\begin{table}
	\centering 
	\begin{tabular}{lr@{\,}r@{\,}lr@{\,}r@{\,}lr@{\,}r@{\,}l} 
		\toprule
		& \multicolumn{9}{c}{Average angular error\textbar Average end point error:  (Seconds)}\\
		\cmidrule{1-4} 
		\cmidrule{5-10}
		& \multicolumn{3}{c}{Dimetrodon}
		& \multicolumn{3}{c}{ Hydrangea}
		& \multicolumn{3}{c}{ Rubberwhale} \\
		\cmidrule{1-10} 
		ALG1&& 2.85\textbar0.15&(12.91s)  && 2.42\textbar0.21&(25.91s) && 4.16\textbar0.13&(23.55s) \\
		PDR && 3.18\textbar0.17&(15.03s)  && 2.07\textbar\underline{0.17}&(33.34s) && \underline{3.13}\textbar\underline{0.10}&(30.35s)\\
		\midrule 
		pADMM &&\underline{2.62}\textbar \underline{0.13}&(8.90s) &&2.08\textbar\underline{0.17}&(25.97s) &&3.46\textbar0.11&(24.12s)\\
		rpADMMI &&2.68\textbar0.14&(8.90s) &&2.07\textbar0.18&(25.68s) &&3.60\textbar0.12&(26.56s)\\
		rpADMMII &&2.71\textbar0.14&(8.90s) &&\underline{2.06}\textbar0.18&(26.01s) &&3.76\textbar0.12&(24.98s)\\
		Zach-pADMM&& 2.96\textbar0.16&(6.76s)  && 2.65\textbar0.21&(13.69s) && 5.11\textbar0.16&(13.57s)\\
		\bottomrule 
	\end{tabular}
	\vspace*{0.5em}
	\caption{ Numerical results for the TV-$L^1$ optical flow estimates.
		In the table, we use $s|l(t)$ with $s$ representing the average angular error (AAE), $l$ denoting the average end point error (EPE), and  $s$ representing the computation time with second. The best results of AAE or EPE are underlined.}
	\label{tab:optical:zachp:beta}
\end{table}

\begin{table}
	\centering
	\label{tab:ani:atq:atv}
	\begin{tabular}{|c|c|c|c|c|c|}
		\hline
		& \multicolumn{5}{c|}{Average angular error\textbar Average end point error}
		\\
		\hline 
		&  Urban2 & Grove2 & Urban3& Venus & Grove3   \\
		\hline
		pADMM & \underline{2.60}\textbar\underline{0.36} & \underline{2.25}\textbar\underline{0.15} & \underline{4.24}\textbar0.54& \underline{4.36}\textbar\underline{0.29} & 6.27\textbar0.66 \\
		rpADMMI & 2.78\textbar\underline{0.36} & \underline{2.25}\textbar\underline{0.15}& 4.25\textbar\underline{0.52}& 4.47\textbar\underline{0.29} & 6.24\textbar0.65 \\	
		rpADMMII  & 2.89\textbar0.37 & 2.28\textbar0.16 & \underline{4.24}\textbar\underline{0.52}& 4.53\textbar0.30 & \underline{6.21}\textbar0.65 \\	
		Zach-pADMM & \underline{2.60}\textbar\underline{0.36}& 2.28\textbar0.16 & 4.87\textbar0.54  & \underline{4.36}\textbar\underline{0.29}  &6.55\textbar\underline{0.64}\\
		Zach-O & 3.06\textbar0.38& 2.31\textbar0.16 & 6.63\textbar0.71  & 5.25\textbar0.35  &6.60\textbar0.72\\
		\hline
	\end{tabular}
	\vspace*{0.5em}
	\caption{Numerical results for the TV-$L^1$ optical flow estimates.
		In the table, we use $s|l$ with $s$ representing the average angular error (AAE) and $l$ denoting the average end point error (EPE). The best results of AAE or EPE are underlined.}
	\label{tab:angular:error1}
\end{table}

\section{Conclusions}\label{sec:conclusion}
We give a systematic studies on the preconditioned ADMM for the TV-$L^1$ optical flow estimates.  We developed several novel and efficient preconditioned  ADMM with convergence guarantee. Various efficient block preconditioners are proposed within the ADMM framework. 
The numerical tests for the TV-$L^1$ optical flow estimates showed that the proposed preconditioned ADMM algorithms have
the potential to bring out appealing benefits and fast algorithms with convergence guarantee.

\noindent

\textbf{Acknowledgements}
H. Sun acknowledges the support of National Natural Science Foundation of China under grant No. \,11701563. He also acknowledges the support of Alexander
von Humboldt Foundation and the support from the program of China Scholarship Council (CSC) under No. 201906365017 during the preparations of this work.

\section{Appendix: Some additional details} \label{sec:appendix}

\begin{proposition} The explicit representations of the resolvent $\bar {\mathfrak{X}}=(I + \sigma \partial F)^{-1}(\tilde {\mathfrak{X}})$ as in \eqref{eq:resolvent}  can be obtained by the Fenchel-Rockafellar duality theory.
\end{proposition}
\begin{proof}
	By the definition of the resolvent, we get
	\begin{equation}\label{eq:argmin:F:resol}
	(\bar w,\bar{\bm d}) = (I + \sigma \partial F)^{-1}(\tilde w,\tilde {\bm d}):=\argmin_{w,\bm {d}} F(w,\bm d) +
	\frac{\|w-\tilde w\|^2}{2\sigma} +  \frac{\|\bm d-\tilde {\bm d}\|^2}{2\sigma}.
	\end{equation}
	Denote $\zeta = [\beta , I_x, I_y]$ and
	\[
	G(w,\bm d) = \frac{\|w-\tilde w\|^2}{2\sigma} +  \frac{\|\bm d-\tilde {\bm d}\|^2}{2\sigma}, \quad K = [\beta I, I_x, I_y].
	\]
	Then the problem in \eqref{eq:argmin:F:resol} can be written as 
	\begin{equation}\label{eq:fenchel:primal}
	min_{{\mathfrak{X}}} L(K {\mathfrak{X}}) + G( {\mathfrak{X}}), \quad L(z) = \lambda \|z+I_t\|_{1}, \quad L(K {\mathfrak{X}}) = F(w,\bm d)=\lambda \|\rho(w,\bm d)\|_{1}.
	\end{equation}
	By the Fenchel-Rockafellar duality \cite{KK}, the problem \eqref{eq:fenchel:primal} is equivalent to 
	the following dual problem
	\begin{equation}\label{eq:fenchel:dual}
	\max_{\Theta}G^*(-K^*\Theta) + L^*(\Theta),
	\end{equation}
	where $L^*(\Theta) = I_{\{\|\Theta\|_{\infty} \leq 1\}} (\Theta) - \langle I_t, \Theta \rangle $ and $G^*(z)=\frac{\sigma}{2}\|z_1\|^2 + \langle \tilde w, z_1 \rangle  + \frac{\sigma}{2}\|z_2\|^2 + \langle \tilde {\bm d}, z_2 \rangle $ through direct calculation.
	The primal 
	solution $\bar {\mathfrak{X}}$ of \eqref{eq:fenchel:primal} and the dual solution $\bar 
	\Theta$ of \eqref{eq:fenchel:dual} have the following 
	optimality conditions 
	\begin{equation}\label{eq:optimality:condi}
	K  \bar {\mathfrak{X}} \in \partial L^*(  \bar \Theta), \quad -K^* \bar \Theta \in \partial G(\bar {\mathfrak{X}}).
	\end{equation}
	With \eqref{eq:fenchel:dual}, we get the optimal solution 
	\begin{equation}
	\bar \Theta = \mathcal{P}_{\Theta}(\frac{\rho(\tilde w, \tilde {\bm d})}{\sigma (\beta^2 + I_x^2 + I_y^2)})= \begin{cases}
	\frac{\rho(\tilde w, \tilde {\bm d})}{\sigma (\beta^2 + I_x^2 + I_y^2)}, \quad |\rho(\tilde w, \tilde{\bm d})| \leq \sigma \lambda (\beta^2 + I_x^2 + I_y^2), \\
	\lambda \rho(\tilde w, \tilde {\bm d})/|\rho(\tilde w, \tilde {\bm d})|,\quad  |\rho(\tilde w, \tilde {\bm d})| > \sigma \lambda (\beta^2 + I_x^2 + I_y^2).
	\end{cases}
	\end{equation}
	where $\mathcal{P}_{\Theta}$ is the projection to the set $\{ \Theta :  \|\Theta\|_{\infty} \leq \lambda\}$ and is the same as $\mathcal{P}_s(\cdot)$. By $ -K^* \bar \Theta \in \partial G(\bar {\mathfrak{X}})$ in \eqref{eq:optimality:condi}, noting that $|\zeta|^2 = \beta^2 + I_x^2 + I_y^2$ and $\partial G(\bar {\mathfrak{X}})$ is single valued, i.e.,   we have $ -K^* \bar \Theta = \nabla G(\bar {\mathfrak{X}})$. Written it component-wisely, we obtain
	\[
	-\begin{bmatrix}
	\beta \\
	I_x \\
	I_y
	\end{bmatrix} \bar \Theta = \begin{bmatrix}
	\frac{\bar w-\tilde w}{\sigma} \\
	\frac{\bar u -\tilde u}{\sigma} \\
	\frac{\bar v -\tilde v}{\sigma} 
	\end{bmatrix} \Rightarrow
	\begin{bmatrix}
	\bar w \\
	\bar u\\
	\bar v
	\end{bmatrix}
	= \begin{bmatrix}
	\tilde w \\
	\tilde u\\
	\tilde v
	\end{bmatrix} +\begin{cases}
	-\rho(\tilde w, \tilde {\bm d}) \zeta^T/|\zeta|^2, \quad \text{if} \ \  |\rho(\tilde w, \tilde {\bm d})| \leq \sigma \lambda (\beta^2 + I_x^2 + I_y^2), \\
	-\lambda \sigma \zeta^T, \quad \text{if} \ \ \rho(\tilde w, \tilde {\bm d}) > \sigma \lambda (\beta^2 + I_x^2 + I_y^2), \\
	\lambda \sigma \zeta^T, \quad \text{if} \ \ \rho(\tilde w, \tilde {\bm d}) <- \sigma \lambda (\beta^2 + I_x^2 + I_y^2). \\
	\end{cases}
	\]
	The proof is finished.
\end{proof}
Now let's turn to the detail of the updates for \ref{eq:rpADMMI:cp}, \ref{eq:rpADMMII:sun}, \ref{eq:admm:zach}. The updates of  \ref{iteration:opflow:l1:pdr} are similar and are thus omitted. For the update of $y^{k+1}$ in \eqref{eq:aug:pre1:vector}, since
\begin{align}
&0 \in \mathcal{B}^* \Lambda^k - (ac I-c\mathcal{B}^*\mathcal{B})(y-y^k) -\partial \mathcal{H}(y) -c \mathcal{B}^*(\mathcal{A}s^k+ \mathcal{B}y) \notag \\
& y^{k+1} =(I + \frac{1}{ac}\partial \mathcal{H})^{-1}(y^k - \frac{1}{a}\mathcal{B}^*\mathcal{B} y^k -\frac{1}{a}\mathcal{B}^*\mathcal{A}s^k +\frac{1}{ac}\mathcal{B}^*\Lambda^k), \label{eq:y}
\end{align}
and
\begin{align*}
&\mathcal{B}^*\mathcal{B} = \text{Diag}[\nabla \nabla^*, \nabla \nabla^*], \quad \mathcal{B}^*\mathcal{A}s^k = [\nabla (\beta s^k), \nabla(I_x s^k), \nabla(I_ys^k)]^{T}, \quad \mathcal{B}^*\Lambda^k = [\nabla u^k, \nabla v^k]^T, \\
&(I + \frac{1}{ac}\partial \mathcal{H})^{-1} = (\mP_p(\cdot),\mP_q(\cdot))^T,
\end{align*}
substituting these to \eqref{eq:y}, we thus get the update of $p^{k+1}$ and $q^{k+1}$ in \eqref{eq:rpADMMI:cp}.  For the update of $s^{k+1}$, since
\begin{align}
&0 \in \mathcal{A}^*\Lambda^k -c\mathcal{A}^*(\mathcal{A}s+\mathcal{B} y^{k+1}) - \partial \mathcal{G}(s) - (\tilde a c I - cM_s)(s-s^k) \notag \\
& s^{k+1} = (I + \frac{1}{\tilde a c}\partial \mathcal{G})^{-1}(s^k - \frac{1}{\tilde a}M_ss^k -\frac{1}{\tilde a}\mathcal{A}^*\mathcal{B} y^{k+1} +\frac{1}{\tilde a c}\mathcal{A}^* \Lambda^k ), \label{eq:s}
\end{align}
and 
\begin{align*}
&-\mathcal{A}^*\mathcal{B} y^{k+1} = \beta \D p^{k+1} + I_x \D q_1^{k+1} + I_y \D q_2^{k+1}, \quad \mathcal{A}^*\Lambda^k = \beta w^k + I_x u^k + I_y v^k, \\
&(I + \frac{1}{\tilde a c}\partial \mathcal{G})^{-1}(\cdot) = \mP_s(\cdot + \frac{1}{\tilde a c}I_t),
\end{align*}
substituting these to equation \eqref{eq:s}, we thus get the update of $s^{k+1}$.

  \end{document}